\documentclass{amsart}

\usepackage{amsthm}
\usepackage{amsfonts}
\usepackage{amsmath}
\usepackage{amssymb}
\usepackage{mathrsfs}
\usepackage{fullpage}
\usepackage{stmaryrd}
\usepackage{hyperref}
\usepackage{verbatim}

\theoremstyle{plain}

\theoremstyle{definition}

\theoremstyle{remark}

\newcommand{\Begin}[2]{\begin{#1}\label{#2}}

\newcommand{\bbQ}{\mathbb{Q}}
\newcommand{\bbR}{\mathbb{R}}

\newcommand{\CA}{\mathcal{A}}

\newcommand{\CM}{\mathcal{M}}
\newcommand{\CN}{\mathcal{N}}

\newcommand{\SCRL}{\mathscr{L}}

\newcommand{\KP}{\mathsf{KP}}

\newcommand{\gbbR}{\bbR^{\geq 0}}
\newcommand{\CD}{\mathcal{D}}

\begin{document}

\title{Bounds on Continuous Scott Rank}

\author{William Chan}
\address{Department of Mathematics, University of North Texas, Denton, TX 76203}
\email{William.Chan@unt.edu}

\author{Ruiyuan Chen}
\address{Department of Mathematics, University of Illinois at Urbana-Champaign, Urbana, IL 61801}
\email{ruiyuan@illinois.edu}

\begin{abstract}
An analog of Nadel's effective bound for the continuous Scott rank of metric structures, developed in \cite{Metric-Scott-Analysis}, will be established: Let $\SCRL$ be a language of continuous logic with code $\hat\SCRL$. Let $\Omega$ be a weak modulus of uniform continuity with code $\hat\Omega$. Let $\CD$ be a countable $\SCRL$-pre-structure. Let $\bar{\CD}$ denote the completion structure of $\CD$. Then $\mathrm{SR}_\Omega(\bar{D}) \leq \omega_1^{\hat\SCRL\oplus\hat\Omega\oplus\CD}$, the Church-Kleene ordinal relative to $\hat\SCRL\oplus\hat\Omega\oplus\CD$.
\end{abstract}

\thanks{July 31, 2019. The first author was supported by NSF grant DMS-1703708.}

\maketitle

\section{Introduction}\label{Introduction}

The authors of \cite{Metric-Scott-Analysis} developed a new Scott analysis for continuous logic for metric structures. If $\SCRL$ is a language of continuous logic, then an $\SCRL$-structure of continuous logic is a Polish metric space endowed with a suitable interpretation for each symbol of $\SCRL$. In this setting, a countable dense $\SCRL$-pre-structure (Definition \ref{pre-structures}) completely determines the original structure through taking completions. The main goal of this paper is to establish a countable effective bound on the continuous Scott rank of a metric structure which depends on the definability content of the countable dense $\SCRL$-pre-structure in a manner analogous to Nadel's effective bound for countable first order structure.

To motivate the ideas of the Scott analysis and the use of continuous logic to study Polish metric structures, a brief review of the classical Scott analysis will be given:

Let $\SCRL$ be a first order language in the classical sense. Let $\CM$ denote some $\SCRL$-structure. Let $\bar{a}$ and $\bar{b}$ be two tuples in $\CM$ of the same length. The Scott analysis begins by attempting to quantify how difficult it is to distinguish $\bar{a}$ and $\bar{b}$ in a manner expressible by $\SCRL$. For instance, if there was a $\SCRL$-automorphism of $\CM$ taking $\bar{a}$ to $\bar{b}$, one would consider $\bar{a}$ and $\bar{b}$ indistinguishable by the first order expressive power of $\SCRL$.

By recursion, one defines for each ordinal $\alpha$, the back-and-forth relation $\sim_\alpha$ on finite tuples of elements of $\CM$ as follows: Let $\bar{a} = (a_0,...,a_{p - 1})$ and $\bar{b} = (b_0,...,b_{p - 1})$ where $p \in \omega$.

\begin{itemize}
\item $\bar{a} \sim_0 \bar{b}$ if and only if the map taking $a_i$ to $b_i$ where $i < p$ is a partial $\SCRL$-isomorphism.

\item Suppose $\sim_\alpha$ has been defined. Then $\bar{a} \sim_{\alpha + 1}\bar{b}$ if and only if $(\forall c)(\exists d)(\bar{a}c\sim_\alpha\bar{b}d) \wedge (\forall d)(\exists c)(\bar{a}c\sim_\alpha\bar{b}d)$. 

\item Suppose $\alpha$ is a limit ordinal and for all $\beta < \alpha$, $\sim_\beta$ has been defined. Then $\bar{a} \sim_\alpha \bar{b}$ if and only if for all $\beta < \alpha$, $\bar{a}\sim_\beta \bar{b}$.
\end{itemize}

Intuitively, if $\neg(\bar{a} \sim_0 \bar{b})$ holds, then one can say that $\bar{a}$ and $\bar{b}$ has been distinguished and in fact, the two tuples fail to satisfy the same atomic formulas. If $\neg(\bar{a} \sim_1 \bar{b})$, then $\bar{a}$ and $\bar{b}$ have been distinguished by failing a property which is expressible by one existential quantification over atomic formulas. For each $\alpha$, $\sim_\alpha$ is closely connected to the type satisfied by tuples via formula of quantifier rank less than or equal to $\alpha$. Of course, if there is an $\SCRL$-automorphism of $\CM$ taking $\bar{a}$ to $\bar{b}$, then $\bar{a} \sim_\alpha \bar{b}$ for all $\alpha \in \mathrm{ON}$. Thus $\bar{a}$ and $\bar{b}$ are indistinguishable by infinitary $\SCRL$-formulas.

Note that $\sim_\alpha$ is an equivalence relation on tuples from $M$. Each tuple is $\sim_\alpha$-related only to tuples of the same length and if $\alpha \leq \beta$, then ${\sim_\beta} \subseteq {\sim_\alpha}$. One definition of the Scott rank states that $\mathrm{SR}(M)$ is the least ordinal $\alpha$ so that for all $\beta \geq \alpha$, ${\sim_\alpha} = {\sim_\beta}$. Intuitively, the Scott rank of $\CM$ is the least ordinal $\alpha$ so that every pair of tuples in $\CM$ which can be distinguished by an infinitary $\SCRL$-formula has been distinguished by an infinitary formula of rank less than or equal to $\alpha$. Each distinguishable pair $(\bar{a},\bar{b})$ corresponds to an infinitary $\SCRL$-formula which makes the distinctions. By roughly collecting all these formulas into a single formula, one obtains a formula $\varphi_\CM$ which is called the Scott sentence for $\CM$. If the language $\SCRL$ is countable and $\CM$ is countable, then this formula is an invariant distinguishing $\CM$ from all other countable $\SCRL$-structure: that is, if $\CN$ is a countable $\SCRL$-structure, then $\CN \models \varphi_\CM$ if and only if $\CM$ and $\CN$ are $\SCRL$-isomorphic.

By a simple cardinality consideration, one has that $\mathrm{SR}(\CM) < |\CM|^+$. One can also view the back-and-forth process as a monotone operator which collects the distinguishable pairs: Let $\Gamma: \mathscr{P}(\CM^{<\omega}\times\CM^{<\omega}) \rightarrow \mathscr{P}(\CM^{<\omega}\times\CM^{<\omega})$ be defined by $\Gamma(X)$ is the set of $(\bar{a},\bar{b}) \in \CM^{<\omega}\times\CM^{<\omega}$ so that one of the following holds:

\begin{enumerate}
\item $|\bar{a}| \neq |\bar{b}|$, that is the length of the two tuples are different.

\item $\neg(\bar{a} \sim_0 \bar{b})$. 

\item $(\exists c)(\forall d)((\bar{a}c,\bar{b}d) \in X) \vee (\exists d)(\forall c)((\bar{a}c,\bar{b}d) \in X)$
\end{enumerate}

By recursion, define $\langle I_\Gamma^\alpha : \alpha \in \mathrm{ON}\rangle$ by: $I_\Gamma^0 = \Gamma(\emptyset)$, $I^{<\alpha}_\Gamma = \bigcup_{\beta < \alpha} I_\Gamma^\beta$, and $I_\Gamma^\alpha = \Gamma(I_\Gamma^{<\alpha})$. The closure ordinal $\|\Gamma_q\|$ of $\Gamma$ is the least ordinal $\alpha$ so that $I^\alpha_\Gamma = I^{<\alpha}_\Gamma$. If $\Gamma$ is a monotone operator which is definable by a positive $\Sigma$-formula in an admissible set $\CA$, then Gandy showed that $\|\Gamma\| \leq o(\CA) = \CA \cap \mathrm{ON}$, the ordinal height of $\CA$. (See \cite{Admissible-Sets-and-Structures} Chapter VI for more on inductive definability in admissible sets.) For the monotone operator $\Gamma$ defined above, one can check that all distinguishable pairs of tuples from the Scott analysis appear in $I_\Gamma^{\|\Gamma\|}$ and $\Gamma$ is positive $\Sigma$-definable in any admissible set containing the language $\SCRL$ and the structure $\CM$. If $\SCRL$ and $\CM$ are both countable (so they can essentially be coded by reals), then the minimal such admissible set is the initial segment of G\"odel's relativized constructible universe, $L_{\omega_1^{\SCRL \oplus \CM}}(\SCRL \oplus \CM)$, which has ordinal height $\omega_1^{\SCRL \oplus \CM}$. $\omega_1^{\SCRL \oplus \CM}$ is the Church-Kleene ordinal relative to $\SCRL\oplus\CM$ which is defined to be the least ordinal that does not have a presentation on $\omega$ recursive in $\SCRL \oplus \CM$. Thus one has obtained Nadel's \cite{Scott-Sentences-and-Admissible-Sets} effective bound which asserts that if $\CM$ is a countable $\SCRL$-structure with $\SCRL$ a countable first order language, then $\mathrm{SR}(\CM) \leq \omega_1^{\SCRL \oplus \CM}$. 

The above is the definition of Scott rank that \cite{Metric-Scott-Analysis} attempts to generalize. There are other variations of the Scott rank that focus in on tuples $\bar{a}$ rather than just the stabilization point of the back-and-forth relations. That is, one can define $\mathrm{SR}^*(\bar{a},\bar{b})$ to be the least ordinal $\alpha$ so that $\neg(\bar{a}\sim_\alpha\bar{b})$ if such an ordinal exists. Otherwise say $\mathrm{SR}^ *(\bar{a},\bar{b}) = \infty$. One can define $\mathrm{SR}^*(\bar{a}) = \sup\{\mathrm{SR}^*(\bar{a},\bar{b}) : \mathrm{SR}^*(\bar{a},\bar{b}) \neq \infty\}$. Then $\mathrm{SR}^*(\CM) = \sup\{\mathrm{SR}^*(\bar{a}) + 1 : \bar{a} \in \CM^{<\omega}\}$. The same argument as above shows that if $\CM$ is a countable structure in a countable language $\SCRL$, then $\mathrm{SR}^*(\CM) \leq \omega_1^{\SCRL \oplus \CM} + 1$. The definition of Scott rank using $\mathrm{SR}^*$ is somewhat more common. It is the form used in \cite{Scott-Rank-of-Polish-Metric-Spaces} and \cite{Bounds-Scott-Ranks-Some-Polish-Metric-Space}. The distinction between recursive structures $\CM$ and $\CN$ in a recursive language $\SCRL$ so that $\mathrm{SR}^*(\CM) = \omega_1^\emptyset + 1$ (for example, the Harrison linear ordering) or $\mathrm{SR}^*(\CN) = \omega_1^\emptyset$ (for example, the linear ordering of Makkai \cite{An-Example-Concerning-Scott-Heights}) are of particular interest in computable model theory. If $\CM$ and $\CN$ are two recursive structures so that $\mathrm{SR}^*(\CM) = \omega_1^\emptyset + 1$ and $\mathrm{SR}^*(\CN) = \omega_1^\emptyset$, then $\mathrm{SR}(\CM) = \mathrm{SR}(\CN) = \omega_1^\emptyset$. The definition of Scott rank used in this article cannot make this distinction; however, for the purpose of finding an effective bound, this will not be relevant.

Polish metric spaces are complete separable metric spaces. The uncountable Polish metric spaces have cardinality $2^{\aleph_0}$. By cardinality considerations, the classical first order Scott rank of a Polish metric space is less than $(2^{\aleph_0})^+$. Although Polish metric spaces are not countable, they are entirely determined by their countable dense metric subspace. A natural question asked by Fokina, Friedman, Koerwien, and Nies \cite{Scott-Rank-Polish-Metric-Spaces-FFKN} was whether the Scott rank of a Polish metric space is countable. The first author in \cite{Bounds-Scott-Ranks-Some-Polish-Metric-Space} asked whether $\mathrm{SR}(\bar{\CM}) \leq \omega_1^{\CM}$, where $\bar{\CM}$ is the completion of the countable metric space $\CM$. This is the natural analog of Nadel's effective bound for Polish metric spaces.

However, the question of Fokina, Friedman, Koerwien, and Nies remains open. (See \cite{Scott-Rank-of-Polish-Metric-Spaces} and \cite{Scott-Rank-of-Polish-Metric-Space-Erratum}.) Some partial results are known. Fokina, Friedman, Koerwien, and Nies \cite{Scott-Rank-Polish-Metric-Spaces-FFKN} showed that if $\CM$ is a compact Polish metric space, then $\mathrm{SR}(\CM) = \omega$ and, in fact, $\mathrm{SR}^*(\CM) = \omega+1$. Nies informed the first author that their argument used some results of Gromov. See \cite{Bounds-Scott-Ranks-Some-Polish-Metric-Space} for a combinatorial proof using the K\"onig's lemma. Doucha \cite{Scott-Rank-of-Polish-Metric-Spaces} showed that if $\CM$ is a Polish metric space, then $\mathrm{SR}^*(\CM) \leq \omega_1$. So the original question becomes whether there exists a Polish metric space $\CM$ of Scott rank exactly $\omega_1$. \cite{Bounds-Scott-Ranks-Some-Polish-Metric-Space} uses admissible sets and infinitary logic in countable admissible fragments to give another proof of Doucha's result and some additional partial results. A metric space is proper if and only if all its closed balls are compact. \cite{Bounds-Scott-Ranks-Some-Polish-Metric-Space} showed that if $\CM$ is a countable metric space so that the completion $\bar{\CM}$ is a proper Polish metric space, then $\mathrm{SR}^*(\bar{\CM}) \leq \omega_1^{\CM} + 1$. If $\CM$ is a countable metric space so that $\bar{\CM}$ is rigid (has no nontrivial autoisometry), then $\mathrm{SR}^*(\CM) < \omega_1^{\CM}$.

The results of \cite{Bounds-Scott-Ranks-Some-Polish-Metric-Space} are proved by using winning strategies in approximation forms of the Ehrenfeucht-Fra\"isse game in an illfounded model of $\KP$. The arguments are quite different than the classical method involving the monotone operator. Although the countable submetric space is essential in representing elements of the completion via Cauchy sequences, the connected between these techniques and the use of the countable dense submetric space of the completion seem very weak. 

An alternative logic that appears more suitable for structures on Polish metric spaces is continuous logic for metric structures. The reader should consult \cite{Model-Theory-for-Metric-Structures} and \cite{Metric-Scott-Analysis} for more details on continuous logic for metric structures. A language $\SCRL$ of continuous logic consists of function, relation, and constant symbols. In addition to the arity, each symbol is associated with a modulus of uniform continuity. Connectives are now certain continuous bounded real-valued functions. The intended structures are Polish metric spaces with the functions and relations interpreted by continuous functions on the Polish metric space respecting the indicated modulus.

The authors of \cite{Metric-Scott-Analysis} proposed and developed a Scott analysis for continuous logic for metric structures. They defined an analogous back-and-forth pseudo-distance which depends on one additional object called a weak modulus of continuity $\Omega$. For each language $\SCRL$, weak modulus $\Omega$, and separable $\SCRL$-structure $\CN$ of continuous logic, they defined a Scott rank $\mathrm{SR}_\Omega(\CN)$. From their Scott analysis, they derived a Scott sentence $\varphi_\CN$ so that for all separable $\SCRL$-structure $\CM$, $\CM \models (\varphi_\CN = 0)$ if and only $\CN$ and $\CM$ are $\SCRL$-isomorphic in continuous logic, in the case that $\Omega$ is a universal weak modulus. See \cite{Metric-Scott-Analysis} Theorem 3.8 and 5.5 for these results and more details. These results give strong evidence that their theory can justly be called a ``Scott analysis'' for continuous logic.

The authors of \cite{Metric-Scott-Analysis} showed that every Polish metric structure in continuous logic has countable Scott rank by cardinality considerations.  Every Polish metric structure in some language $\SCRL$ is the completion $\SCRL$-structure of a countable dense $\SCRL$-pre-structure. The main task addressed in this article is to investigate the connection between the Scott rank of a Polish metric structure and any of its countable dense pre-structures. 

Let $\SCRL$ be a countable $\SCRL$-structure. Let $\hat{\SCRL}$ denote a real coding $\SCRL$ which includes information about the associated modulus of the language. Let $\Omega$ be a weak modulus of continuity which is coded by a real $\hat{\Omega}$. Let $\mathcal{D}$ denote a countable $\SCRL$-pre-structure. The main theorem is
\\*
\\*\noindent\textbf{Theorem \ref{bounds on continuous scott rank}.}\textit{
Let $\SCRL$ be a countable language with code $\hat\SCRL$. Let $\Omega$ be a weak modulus with code $\hat\Omega$. Let $\CD$ be a countable $\SCRL$-pre-structure. Let $\bar{\CD}$ denote the completion structure of $\CD$. Then $\mathrm{SR}_\Omega(\bar\CD) \leq \omega_1^{\hat\SCRL\oplus\hat\Omega\oplus\CD}$. 
}
\\*
\\*\indent As an example, consider the class of metric structures consisting of the pure Polish metric spaces. That is, the language is $\SCRL = \emptyset$, which consists of no additional non-logical symbols. Thus $\SCRL$ has a code which is recursive. One can take $\Omega$ to be the universal weak modulus for this language, which exists and has a recursive code from inspecting the proof of \cite{Metric-Scott-Analysis} Proposition 5.3. Now suppose $\CD$ is a countable metric space. $\CD$ is naturally a pure metric pre-structure of continuous logic. Let $\bar{\CD}$ be the completion of $\CD$ as a pure Polish metric structure. Then Theorem \ref{bounds on continuous scott rank} states the continuous Scott rank of $\bar{\CD}$, $\mathrm{SR}_\Omega(\bar{\CD})$, is less than or equal to $\omega_1^{\CD}$, the Church-Kleene ordinal relative to $\CD$. This establishes the analog of Nadel's effective bound on the continuous Scott rank of Polish metric spaces.

The basic template for the proof of Theorem \ref{bounds on continuous scott rank} is the same as the classical first order argument. Now one attempts to define various monotone operators on $\mathcal{D}^{<\omega} \times \mathcal{D}^{<\omega}$ that are positive $\Sigma$-definable in an admissible set containing $\mathcal{D}$, $\hat{\Omega}$, and $\hat{\SCRL}$. However, for this to be meaningful, one needs to ensure the computation of the back-and-forth pseudo-distance within the desired admissible set evaluates to the correct or true computation as performed in the real world. This amounts to showing that the first back-and-forth function $r_0$ is computed correctly by the appropriate admissible set. This will be shown by producing a countable collection of basic formulas respecting the weak modulus $\Omega$ which has a code in every admissible set containing $\hat{\SCRL}\oplus\hat{\Omega}\oplus\mathcal{D}$ and such that this collection is dense in the collection of all basic formulas respecting the weak modulus under the uniform norm. 

The authors would like to thank Andr\'e Nies for commenting on an early draft of this paper.

\section{Continuous Logic}
See \cite{Model-Theory-for-Metric-Structures} for a more detailed exposition on continuous logic for metric structures.

\Begin{definition}{modulus}
\cite{Metric-Scott-Analysis} A modulus of arity $n$ is a continuous function $\Delta : (\gbbR)^n \rightarrow \gbbR$ such that $\Delta(\bar{0}) = 0$, and for all $\bar r,\bar s \in (\gbbR)^n$, $\Delta(\bar r) \leq \Delta(\bar r + \bar s) \leq \Delta(\bar r) + \Delta(\bar s)$. ($+$ refers to coordinate-wise addition.) The latter succinctly states that $\Delta$ is non-decreasing and subadditive.

Let $\Delta$ be a modulus of arity $n$. Let $(X_i,d_{X_i})$ and $(Z,d_Z)$ be metric spaces where $i < n$. Define $d^\Delta : (\prod_{i < n} X_i)^2 \rightarrow \gbbR$ by
$$d^\Delta(\bar{x},\bar{y}) = \Delta((d_{X_0}(\bar{x}(0), \bar{y}(0)), ... , d_{X_{n - 1}}(\bar{x}(n - 1), \bar{y}(n - 1))))$$
Let $f : \prod_{i < n} X_i \rightarrow Z$. $f$ respects $\Delta$ if and only if $d_Z(f(\bar{x}),f(\bar{y})) \leq d^\Delta(\bar{x},\bar{y})$. 
\end{definition}

\Begin{definition}{weak modulus}
\cite{Metric-Scott-Analysis} A weak modulus is a function $\Omega : (\gbbR)^\omega \rightarrow [0,\infty]$ which is non-decreasing, subadditive, lower semi-continuous in the product topology, separately continuous in each coordinate, and $\Omega(\bar{0}) = 0$.

For each $n \in \omega$, define $\Omega_n : (\gbbR)^n \rightarrow \gbbR$ by $\Omega_n(x_0, ..., x_{n - 1}) = \Omega(x_0, ..., x_{n - 1}, 0,0,0...)$. 
\end{definition}

\Begin{fact}{weak modulus and partial cuts of weak modulus}
(\cite{Metric-Scott-Analysis} Lemma 2.3.) Let $\Omega$ be a weak modulus. For all $n \in \omega$, $\Omega_n$ is a modulus of arity $n$. For all $\bar{r} \in (\gbbR)^\omega$, $\Omega(\bar{r}) = \sup_{n \in \omega} \Omega_n(\bar{r} \upharpoonright n)$. 
\end{fact}

\Begin{definition}{respecting weak modulus}
A function $f : X^n \rightarrow \bbR$ respects the weak modulus $\Omega$ if and only if $f$ respects the modulus $\Omega_n$. 
\end{definition}

\Begin{definition}{language of continuous logic}
A relation symbol consists of a symbol $R$, a natural number $a(R)$, and a modulus $\Delta_R$ of arity $a(R)$. A function symbol consists of a symbol $f$, a natural number $a(f)$, and a modulus $\Delta_f$ of arity $a(f)$. 

A language of continuous logic is a collection $\SCRL$ of relation, function, and constant symbols along with a distinguished binary relation symbol $d$, which is intended to represent the distance function.  

For convenience, one will assume all connectives, relations symbols, and the distance relation can only be interpreted to take value in the interval $[0,1]$.
\end{definition}

\Begin{definition}{pre-structures}
Let $\SCRL$ be a language of continuous logic. A $\SCRL$-pre-structure is a collection $\CM$ consisting of the following: There is a (possibly incomplete) metric space $M$. For each relation symbol $R \in \SCRL$, there is a continuous function $R^\CM : M^{a(R)} \rightarrow [0,1]$ which respects the modulus $\Delta_R$. For each function symbol $f \in \SCRL$, there is a continuous function $f : M^{a(f)} \rightarrow M$ respecting the modulus $\Delta_f$. 

A $\SCRL$-structure is an $\SCRL$-pre-structure where $M$ is a complete metric space.
\end{definition}

\Begin{fact}{completion of a pre-structure}
If $\CD$ is a $\SCRL$-pre-structure, then there is a canonical $\SCRL$-structure on $\bar{D}$, the completion of $D$, obtained by extending all the interpretation of symbols to the completion. This structure is denoted $\bar\CD$. 
\end{fact}

\Begin{definition}{formula continuous logic}
Let $\SCRL$ be a countable language of continuous logic. Fix an infinite set $\langle v_i : i \in \omega\rangle$ of variables. 

The collection of $\SCRL$-terms is the smallest set closed under the following:

\noindent 1. Each $v_i$ is a term. 

\noindent 2. If $f$ is a function symbol and $t_0,...,t_{a(f) - 1}$ are terms, then $f(t_0, ..., t_{a(f) - 1})$ is a term. 

The atomic formulas are generated in the following way: If $R$ is a relation symbol and $t_0, ..., t_{a(R) - 1}$ are terms, then $R(t_0, ..., t_{a(R) - 1})$ is an atomic formula. If $t_1$ and $t_2$ are terms, then $d(t_1,t_2)$ is an atomic formula. 

The collection of $\SCRL$-formulas, denoted $\SCRL_{\omega,\omega}$, is the smallest collection closed under the following: 

\noindent 1. All atomic formulas are formulas.

\noindent 2. If $u : [0,1]^n \rightarrow [0,1]$ is a continuous function and $\varphi_0, ..., \varphi_{n - 1}$ are formulas, then $u(\varphi_0, ..., \varphi_{n - 1})$ is a formula. 

\noindent 3. If $\varphi$ is a formula and $v_i$ is a variable, then $\sup_{v_i} \varphi$ and $\inf_{v_i} \varphi$ are formulas. 
\end{definition}

Each term or formula has a canonically associated modulus. (See \cite{Metric-Scott-Analysis} Section 2.2 for more details.)

\Begin{definition}{basic formulas}
The collection of basic formulas is the smallest collection of formulas closed under 1 and 2 in the definition of $\SCRL_{\omega,\omega}$. 
\end{definition}

\Begin{fact}{standard form of basic formula}
If $\varphi$ is a basic formula, then there are atomic formulas $\varphi_0, ..., \varphi_{k - 1}$ and a continuous function $u : [0,1]^k \rightarrow [0,1]$ so that $\varphi = u(\varphi_0, ..., \varphi_{k - 1})$. 
\end{fact}

\begin{proof}
This is proved by induction.
\end{proof}

\Begin{definition}{interpretation of formula}
Let $\SCRL$ be a language of continuous logic. Let $\CD$ be a $\SCRL$-pre-structure. 

The interpretation of the terms of $\SCRL$ are defined as follows:

\noindent 1. For each variable $v_i$, $v_i^\CD(a_0, ..., a_{n - 1}) = a_i$, if $i < n$ and $a_0, ..., a_{n - 1} \in D$. 

\noindent 2. Suppose $t_0, ..., t_{n - 1}$ are terms mentioning variables $v_0, ..., v_{k - 1}$, $f$ is a $n$-ary function symbol, $\bar{a} = (a_0, ..., a_{k - 1})$ is a tuple from $D$, and each $t_j^\CD(\bar{a})$ has already been defined, then $(f(t_0, ..., t_{n - 1}))^\CD(\bar{a}) = f^\CD(t_0^\CD(\bar{a}), ..., t_{n - 1}^{\CD}(\bar{a}))$. 

The interpretation of formulas is defined recursively as follows:

\noindent 1. If $R$ is a $n$-ary function symbol and $t_0, ..., t_{n - 1}$ are terms, then $(R(t_0, ..., t_{n - 1}))^\CD(\bar{a}) = R^\CD(t_0^\CD(\bar{a}), ..., t_{n - 1}^\CD(\bar{a}))$. 

\noindent 2. If $u : [0,1]^n \rightarrow [0,1]$ is a continuous function and $\varphi_0, ..., \varphi_{n - 1}$ are formulas such that $\varphi_i^\CD(\bar{a})$ has been defined, then $(u(\varphi_0, ..., \varphi_{n - 1}))^\CD(\bar{a}) = u(\varphi_0^\CD(\bar{a}), ..., \varphi^\CD_{n - 1}(\bar{a}))$. 

\noindent 3. Suppose $\varphi$ is a formula such that $\varphi^\CD(\bar{a})$ has been defined, then
$$(\sup_{v_i} \varphi)^\CD(a_0, ..., a_{k - 1}) = \sup_{x \in D}\varphi^\CD(a_0, ..., a_{i - 1}, x, a_{i + 1}, ..., a_{k - 1})$$
\end{definition}

\Begin{definition}{coding definition}
A continuous function $f : (\gbbR)^n \rightarrow \gbbR$ is coded by a function $\hat{f} : \bbQ^n \rightarrow \bbQ^\omega$ with the property that for all $\bar{p} \in (\bbQ^{\geq 0})^n$, $\hat{f}(\bar{p})$ is a Cauchy sequence representing $f(\bar{p})$. Similar coding can be defined if $f : [0,1]^n \rightarrow [0,1]$. 

Let $\SCRL$ be a countable language for continuous logic. A code $\hat{\SCRL}$ for $\SCRL$ consists of the following objects:

\noindent 1. The symbols of $\SCRL$.

\noindent 2. For each function symbol or relation symbol $P \in \SCRL$, $\hat{\Delta}_P$. 

If $\Omega$ is a weak modulus, the code for $\Omega$ is $\hat{\Omega}$ consisting of $\hat{\Omega}_n$ for each $n \in \omega$. 

Each connective $u : [0,1]^n \rightarrow [0,1]$ has a code $\hat{u}$ as defined above. Using this and the recursive definition of formulas, one can obtain codes $\hat{\varphi}$ for each formula $\varphi$.

Let $\CD = (D,d)$ (where $d$ refers to the metric) be a countable $\SCRL$-pre-structure. If $f : D^n \rightarrow \bbR$, then the code for $f$ is $\hat{f} : D^n \rightarrow \bbQ^\omega$ so that for all $\bar{x} \in D^n$, $\hat{f}(\bar{x})$ is a Cauchy-sequence representing $f(\bar{x})$. 

The code of $\CD$, denoted $\hat\CD$, consists of the underlying set $D$, $f^\CD$, $\hat{R}^\CD$, and $\hat{d}$.

If $E$ is a countable collection of functions (which have a code as above), then a code for $E$ is a function $\hat{E}$ on $\omega$ so that for all $n$, $\hat{E}(n)$ is a code for a function in $E$ and for every function $f \in E$, there is some $n$ so that $\hat{E}(n)$ is a code for $f$. 
\end{definition}

\Begin{remark}{coding definition remark}
Throughout, some of coding details will be left to the reader. For instance, the reader can check that in models of $\KP$, values of functions on appropriate objects can be recovered from the appropriate codes of the function and the objects. 

In some sense, information about the $\SCRL$-structure $\bar\CD$ is entirely contained in $\hat\SCRL \oplus \CD$. In the following, one will be concerned about analyzing the continuous Scott analysis of $\bar\CD$. Unlike the first order case, the continuous Scott analysis has one additional parameter, a weak modulus $\Omega$. Hence the real $\hat{\SCRL}\oplus\hat{\Omega}\oplus\CD$ codes all the parameters in the continuous Scott analysis of $\CM$ relative to the weak modulus $\Omega$. The main concern is to find a bound on the continuous Scott rank of $\CM$ relative to $\Omega$ using these parameters. 
\end{remark}

\Begin{definition}{formulas of some admissible set}
Let $\SCRL$ be a countable language. Let $\Omega$ be a weak modulus. Let $\CD$ be a countable $\SCRL$-pre-structure (which one may assume has domain $\omega$). Let $\CA$ be a countable admissible set. 

$\SCRL$ (respectively, $\Omega$) is said to belong to the admissible set $\CA$ if and only if $\hat\SCRL \in A$ (or $\hat{\Omega} \in A$, respectively).

$\CD$ is said to belong to $\CA$ if and only if $\CD \in A$. 

If $\varphi \in L_{\omega,\omega}$, then $\varphi$ belongs to $\CA$ if and only if $\hat{\varphi} \in A$. 

Let $\SCRL_{\omega,\omega}^\CA$ denote the collection of all formulas of $\SCRL_{\omega,\omega}$ that belong to $\CA$. 
\end{definition}

The interpretation of $\SCRL_{\omega,\omega}$ formulas evaluated at elements of $\CD$ is the same as the interpretation in the completion $\bar\CD$:

\Begin{fact}{interpretation prestructure agree on completion}
Let $\CD$ be a countable $\SCRL$-pre-structure. Let $x_0, ..., x_{n - 1} \in D$. Let $\varphi \in \SCRL_{\omega,\omega}$. Then $\varphi^\CD(x_0, ..., x_{n - 1}) = \varphi^{\bar\CD}(x_0, ..., x_{n - 1})$. 
\end{fact}

\begin{proof}
For all terms $t$ and $\bar{a} \in D$, $t^\CD(\bar{a}) = t^{\bar{D}}(\bar{a})$. 

Fact \ref{interpretation prestructure agree on completion} is proved by induction:

Suppose the result holds for $\varphi_0, ..., \varphi_{n - 1}$ and $u : [0,1]^n \rightarrow [0,1]$ is a continuous function, then
$$(u(\varphi_0, ..., \varphi_{ n- 1}))^{\bar\CD}(\bar{a}) = u(\varphi^{\bar\CD}_0(\bar{a}), ..., \varphi_{n - 1}^{\bar{\CD}}(\bar{a})) = u(\varphi_0^\CD(\bar{a}), ..., \varphi_{n - 1}^\CD(\bar{a})) = (u(\varphi_0, ..., \varphi_{n - 1}))^\CD(\bar{a})$$
by using the induction hypothesis.

Suppose the result holds for $\varphi$. Let $i < k$. Then
$$(\sup_{v_i}\varphi)^{\bar\CD}(a_0, ..., a_{k - 1}) = \sup_{b \in \bar{D}} \varphi^{\bar\CD}(a_0, ..., a_{i - 1}, b, a_{i + 1} ..., a_k)$$
By the continuity of $\varphi$ and the induction hypothesis, 
$$=\sup_{b \in D} \varphi^{\bar\CD}(a_1, ..., a_{i - 1}, b, a_{i + 1}, ..., a_{k - 1}) = \sup_{b \in D}\varphi^\CD(a_1, ..., a_{i - 1},b,a_{i + 1}, ..., a_{k - 1}) = (\sup_{v_i}\varphi)^\CD(a_0, ..., a_{k - 1}).$$
A similar argument holds for $\inf$. 

By induction, the result has been shown.
\end{proof}

$\mathsf{KP}$ is capable of formulating the syntax and semantics of continuous logic using the codes of various continuous functions and Cauchy sequence representations of reals and elements of $\bar{D}$. The following is straightforward coding:

\Begin{fact}{admissible has correct interpretation}
Let $\SCRL$ be a countable language of continuous language. Let $\CA$ be a countable admissible set such that $\SCRL \in \CA$. Let $\CD$ be a countable $\SCRL$-pre-structure in $\CA$. Let $\varphi \in L_{\omega,\omega}^\CA$. Let $\bar{a}$ be a tuple of elements in $D$. Let $(\varphi^\CD(\bar{a}))^\CA$ be the computation of $\varphi^\CD(\bar{a})$ (considered as a $\bbQ$-Cauchy sequence) in the admissible set $\CA$. Then $\varphi^\CD(\bar{a}) = (\varphi^\CD(\bar{a}))^\CA$ and by Fact \ref{interpretation prestructure agree on completion}, this is equal to $\varphi^{\bar{\CD}}(\bar{a})$. 
\end{fact}

\section{Density of Basic Formulas Respecting a Weak Modulus}\label{Density of Basic Formulas Respecting a Weak Modulus}

The Scott analysis for continuous logic of \cite{Metric-Scott-Analysis} is formulated using the back-and-forth pseudo-distance functions $\langle r_\alpha : \alpha \in \mathrm{ON}\rangle$ from Definition \ref{back-and-forth pseudo-distance}. Like the classical Scott analysis, the continuous Scott analysis can be expressed through a certain monotone operator. The Scott rank is then the closure ordinal of this operator. To obtain an effective bound, one will show that there is an equivalent monotone operator in a suitable countable admissible set that correctly represents the true Scott analysis occurring in the real world.

This amounts to showing that one can compute each $r_\alpha$ correctly in an appropriate admissible set. From Definition \ref{back-and-forth pseudo-distance}, one can see that from $r_\alpha$, one obtains $r_{\alpha + 1}$ by a simple explicit procedure. Similarly, if $\beta$ is a limit ordinal, one can obtained $r_\beta$ by a simple explicit procedure from the collection $\{r_\gamma : \gamma < \beta\}$. 

Thus one will need to show that there is an appropriate countable admissible set which can compute the initial function $r_0$ correctly. From Definition \ref{back-and-forth pseudo-distance}, $r_0$ is defined by taking a certain supremum over all basic $\SCRL$-formulas respecting a modulus $\Omega$. The collection of all such formulas is uncountable and so certainly does not belong to any countable admissible set. The main result of this section is to show that there is a countable collection of such formulas which is dense in the uniform norm among the collection of all basic $\SCRL$-formulas respecting $\Omega$ and this collection is arithmetical in the code for $\SCRL$ and $\Omega$. See Fact \ref{dense set of Omega respecting basic formulas} for the precise statement.

\Begin{definition}{nice subset gbbR}
Suppose $n \in \omega$. A subset $A \subseteq (\gbbR)^n$ is said to be nice if and only if 

\noindent 1. For all $\bar{x} \in A$ and all $\bar{y} \in (\gbbR)^n$ such that for all $i < n$, $\bar{y}(i) \leq \bar{x}(i)$, then $\bar{y} \in A$. 

\noindent 2. There is an open set $U \subseteq (\gbbR)^n$ with $\bar{0} \in U$ and $U \subseteq A$. 
\end{definition}

The following will be useful notation throughout:

\Begin{definition}{pi function}
For $\bar{x} \in \bbR^k$, let $\pi(\bar{x}) \in (\gbbR)^k$ be defined by $\pi(\bar{x})(i) = |\bar{x}(i)|$ for all $i < k$.
\end{definition}

\Begin{fact}{largest modulus below a function}
Suppose $f : A \rightarrow \gbbR$ is a function defined on a nice $A \subseteq (\gbbR)^n$, is non-decreasing, $f(\bar{0}) = 0$, continuous at $\bar{0}$, and has code $\hat{f}$. Then there is a function $g : (\gbbR)^n \rightarrow \gbbR$ so that

(i) $g$ is a modulus of arity $n$ with the property that for all $\bar{x} \in A$, $g(\bar{x}) \leq f(\bar{x})$. 

(ii) $g$ is the largest modulus of arity $n$ below $f$ in the following sense: if $h$ is a modulus of arity $n$ below $f$ in the sense that for all $\bar{x} \in A$, $h(\bar{x}) \leq f(\bar{x})$, then one has that for all $\bar{x} \in (\gbbR)^n$, $h(\bar{x}) \leq g(\bar{x})$.

(iii) $\hat{g}$, the code of $g$, is arithmetic in any function $\underline{f} : A \cap (\bbQ^{\geq 0})^n \rightarrow \bbQ^\omega$ such that for all $\bar{p} \in A \cap (\bbQ^{\geq 0})^n$, $\underline{f}(\bar{p})$ is a Cauchy sequence representing $f(\bar{p})$. 

(Note that $\underline{f}$ is essentially a real and it codes $f$ on $(\bbQ^{\geq0})^n$ in a manner similar to how continuous functions are coded, but if $f$ is not continuous, then $f$ may not be recovered from $\underline{f}$.)
\end{fact}

\begin{proof}
For notational simplicity, if $\bar{x}$ and $\bar{y}$ are elements of $(\gbbR)^n$, then one writes $\bar{x} \preceq \bar{y}$ if and only if for all $i < n$, $\bar{x}(i) \leq \bar{y}(i)$. 

Define function $g$ on $(\gbbR)^n$ as follows:
$$g(\bar{x}) = \inf\left\{f(\bar x_0) + ... + f(\bar x_{k - 1}) : k \in \omega \wedge \bar x_0, ..., \bar x_{k - 1} \in A \cap (\bbQ^{\geq 0})^n \wedge \bar x \preceq \sum_{i < k} \bar x_i\right\}$$
Note that $g(\bar 0) = 0$ since $f(\bar{0}) = 0$ and for all $\bar x \in A \cap (\bbQ^{\geq 0})^n$, $g(\bar x) \leq f(\bar x)$.

(Claim 1) $g$ is non-decreasing: Suppose $\bar{x},\bar{y} \in (\gbbR)^n$ is such that $\bar{x} \preceq \bar{y}$. Fix some $\bar{y}_0, ..., \bar{y}_{k - 1} \in A \cap (\bbQ^{\geq 0})^n$ so that $\bar{y} \preceq \sum_{i < k} \bar{y}_i$. Then $\bar{x} \preceq \sum_{i < k}\bar{y}_k$ since $\bar{x} \preceq \bar{y}$. Since $\bar{y}_0,...,\bar{y}_{k - 1}$ with $\bar{y} \leq \sum_{i < k}\bar{y}_i$ were arbitrary, one has that $g(\bar{x}) \leq g(\bar{y})$ by the definition of $g$. 

(Claim 2) $g$ is subadditive: Fix $\bar{x}$ and $\bar{y}$. Suppose $\bar{x}_0, ..., \bar{x}_{k - 1} \in A \cap (\bbQ^{\geq 0})^n$ are such that $\bar{x} \preceq \sum_{i < k} \bar{x}_i$. Suppose $\bar{y}_0, ..., \bar{y}_{p - 1} \in A \cap (\bbQ^{\geq 0})^n$ are such that $\bar{y} \preceq \sum_{i < p}\bar{y}_i$. Let $\bar{z}_0, ..., \bar{z}_{k + p - 1}$ enumerate $\bar{x}_0, ..., \bar{x}_{k - 1}, \bar{y}_0, ..., \bar{y}_{p - 1}$. Then $\bar{x} + \bar{y} \preceq \sum_{i < k + p} \bar{z}_i$. It follows from the definition of $g$ that $g(\bar{x} + \bar{y}) \leq g(\bar{x}) + g(\bar{y})$. 

(Claim 3) $g$ is continuous: Fix $\epsilon > 0$. Since $f$ is continuous at $\bar{0}$ and property 2 of the niceness of $A$, there is some $\delta > 0$ so that for all $\bar x \in (\gbbR)^n$ with the property that for all $i < n$, $\bar x(i) < \delta$, one has that $\bar{x} \in A$ and $f(\bar x) < \epsilon$. 

Fix $\bar{x} \in (\gbbR)^n$. Let $y \in (\gbbR)^n$ be such that for all $i < n$, $\pi(\bar{y} - \bar{x})(i) < \frac{\delta}{2}$. Choose a $\bar{z} \in (\bbQ^{\geq 0})^n$ such that $\pi(\bar{y} - \bar{x}) \preceq \bar{z}$ and $\bar{z}(i) < \delta$ for all $i < k$. By the choice of $\delta$, $\bar{z} \in A \cap (\bbQ^{\geq 0})^n$. Without loss of generality, suppose that $g(\bar{y}) \geq g(\bar{x})$. Note that $\bar{y} \preceq \bar{x} + \pi(\bar{y} - \bar{x})$. Using the fact that $g$ is non-decreasing, subadditive, $g$ is less than $f$ on $A \cap (\bbQ^{\geq 0})^n$, and the choice of $\bar{z}$ and $\delta$, one has that
$$|g(\bar{y}) - g(\bar{x})| = g(\bar{y}) - g(\bar{x}) \leq g(\bar{x} + \pi(\bar{y} - \bar{x})) - g(\bar{x}) \leq g(\bar{x}) + g(\pi(\bar{y} - \bar{x})) - g(\bar{x})$$
$$= g(\pi(\bar{y} - \bar{x})) \leq g(\bar{z}) \leq f(\bar{z}) < \epsilon$$
The continuity of $g$ has been established

It has been shown that $g$ is a modulus. 

(Claim 4) For all $\bar{x} \in A$, $g(\bar{x}) \leq f(\bar{x})$: For each $j \in \omega$, let $\bar{x}_j \preceq \bar{x}$ be such that $\bar{x}_j \in (\bbQ^{\geq 0})^n$ and $(\bar{x} - \bar{x}_j)(i) < \frac{1}{j + 1}$ for all $i < n$. Since $\bar{x} \in A$ and $A$ is nice, $\bar{x}_j \in A \cap (\bbQ^{\geq 0})^n$ for all $j \in \omega$. Since $g(\bar{y}) \leq f(\bar{y})$ for all $\bar{y} \in A \cap (\bbQ^{\geq 0})^n$, the continuity of $g$ and the fact that $f$ is non-decreasing imply that $g(\bar{x}) = \lim_{j \rightarrow \infty} g(\bar{x}_j) = \lim_{j \rightarrow \infty} f(\bar{x}_j) \leq f(\bar{x})$.

(Claim 5) $g$ is the largest $n$-ary modulus below $f$: Suppose $h$ is modulus below $f$ but there is some $\bar{x} \in (\gbbR)^n$ so that $g(\bar{x}) < h(\bar{x})$. Then there is some $\bar{x}_0, ..., \bar{x}_{k - 1} \in A \cap (\bbQ^{\geq 0})^n$ so that $\bar{x} \preceq \sum_{i < k} \bar{x}_i$ and $\sum_{i < k} f(\bar{x}_i) < h(\bar{x})$. However since $h$ is an $n$-ary modulus, one must have
$$\sum_{i < k} f(\bar{x}_i) < h(\bar{x}) \leq \sum_{i < k} h(\bar{x}_i) \leq \sum_{i < k } f(\bar{x}_i)$$
Contradiction.

Note that the definition of $g$ depends only the value of $f$ on $A \cap (\bbQ^{\geq 0})^n$. This implies that there is a code $\hat{g}$ of $g$ which is arithmetic in $\underline{f}$.
\end{proof}

\Begin{fact}{largest modulus to respect a modulus}
Let $\varphi_0, ..., \varphi_{k - 1}$ be atomic formulas of $\SCRL_{\omega,\omega}$ with free variables $v_0, ..., v_{n - 1}$. Let $\Delta_0, ..., \Delta_{k - 1}$ be the canonical modulus of $\varphi_0 ,..., \varphi_{k - 1}$, respectively, as mentioned at the end of Definition \ref{formula continuous logic}. Assume that each $\Delta_i$ is not the constant $0$ function. Let $\Omega$ be a weak modulus. Then there is an $k$-ary modulus $\Delta$ with a code which is arithmetic in $\hat\SCRL\oplus\hat\Omega$ such that for all continuous functions $u : [0,1]^k \rightarrow [0,1]$, $u(\varphi_0, ..., \varphi_k)$ respects $\Omega$ if and only if $u$ respects $\Delta$. 
\end{fact}

\begin{proof}
First, observe that since each $\varphi_i$ is atomic, $\Delta_i$ is recursive in $\hat\SCRL$. Suppose $u : [0,1]^k \rightarrow [0,1]$ is continuous with canonical modulus $\Delta^u$. The canonical modulus for $u(\varphi_0, ..., \varphi_{k - 1})$ is $\Delta^u(\Delta_0, ..., \Delta_{k - 1})$. Suppose $u(\varphi_0, ..., \varphi_{k - 1})$ respects $\Omega_n$. This means that for all $\bar{x} \in (\gbbR)^n$,
$$\Delta^u(\Delta_0(\bar{x}), ..., \Delta_{k - 1}(\bar{x})) \leq \Omega_n(\bar{x})$$ 
Hence $\Delta^u$ must satisfy the following relation: For all $\bar{r} \in (\gbbR)^k$ such that there is some $\bar{x}$ with the property that for all $i < k$, $\Delta_i(\bar{x}) \geq \bar{r}(i)$, 
$$\Delta^u(\bar{r}) \leq \inf\{\Omega_n(\bar{x}) : \bar{x} \in (\gbbR)^n \wedge (\forall i  < k)(\Delta_i(\bar{x}) \geq \bar{r}(i))\}$$

Let $f$ be defined by
$$f(\bar{r}) = \inf\{\Omega_n(\bar{x}) : \bar{x} \in (\gbbR)^n \wedge (\forall i  < k)(\Delta_i(\bar{x}) \geq \bar{r}(i))\}.$$
Using the assumption that each $\Delta_i$ is not constantly $0$, $f$ is defined on a nice $A \subseteq (\gbbR)^k$.

Note that $\Delta^u(\bar{r}) \leq f(\bar{r})$ for all $\bar{r} \in A$ if and only if $u(\varphi_0, ..., \varphi_{k - 1})$ respect $\Omega_n$: 

$(\Rightarrow)$ is clear. 

$(\Leftarrow)$ Suppose there is some $\bar{r} \in A$ so that $f(\bar{r}) < \Delta^u(\bar{r})$. Then there is some $\bar{x}$ so that for all $i$, $\Delta_i(\bar{x}) \geq \bar{r}(i)$ and $\Omega_n(\bar{x}) < \Delta^u(\bar{r})$. Let $\bar{s}$ be such that $\bar{s}(i) = \Delta_i(\bar{x})$. Since for all $i < k$, $\bar{r}(i) \leq \bar{s}(i)$, one has that $\Delta^u(\bar{r}) \leq \Delta^u(\bar{s})$. But then $\Omega_n(\bar{x}) < \Delta^u(\bar{r}) \leq \Delta^u(\bar{s}) = \Delta^u(\Delta_0(\bar{x}), ..., \Delta_{k - 1}(\bar{x}))$. The canonical modulus for the formula does not respect $\Omega_n$. 

$f$ is clearly non-decreasing where it is defined.

Next to show $f$ is continuous at $\bar{0}$: Let $\epsilon > 0$. Since $\Omega_n$ is continuous, there is some $\delta > 0$ so that $\Omega_n(\bar{x}) < \epsilon$ whenever $\bar{x}$ has the property that for all $i < k$, $\bar{x}(i) < \delta$. Since each $\Delta_j$ is not constant in a neighborhood of $\bar 0$, for each $j$, there is some $\bar{z}_j$ with $\bar{z}_j(i) < \frac{\delta}{n}$ for each $i < k$ and $\Delta_j(\bar{z}_j) > 0$. Let $\bar{z} = \sum_{j < k} \bar{z}_j$. Since each $\Delta_i$ is non-decreasing for $0 \leq i < k$, $\Delta_i(\bar{z}) > 0$. Let $\gamma = \min\{\Delta_i(\bar{z}) : i < k\}$. Suppose $\bar{r}$ is such that for all $i < k$, $\bar{r}(i) < \gamma$. Then 
$$\Omega_n(\bar{z}) \in \{\Omega_n(\bar{x}) : \bar{x} \in (\gbbR)^n \wedge (\forall i < k)(\Delta_i(\bar{x}) \geq \bar{r}(i))\}$$
Hence $f(\bar{r}) \leq \Omega_n(\bar{z}) < \epsilon$. $f$ is continuous at $\bar{0}$. 

One can find a function $\underline{f} : A \cap (\bbQ^{\geq 0})^n \rightarrow \bbQ^\omega$ with the properties in Fact \ref{largest modulus below a function} for this function $f$ which is arithmetic in $\hat\SCRL\oplus\hat\Omega$. (Note that in the definition of $f$, one obtains the same function if the infimum is taken over $\bar{x} \in (\bbQ^{\geq 0})^n$ with the required property above.) Fact \ref{largest modulus below a function} states there is a largest modulus $\Delta$ below $f$ which is arithmetic in $\hat\SCRL\oplus\hat\Omega$. This completes the proof.
\end{proof}

\Begin{definition}{lattice}
Let $I \subseteq \bbR$ be an open or closed interval. Let $X$ be a compact metric space. Let $C(X,I)$ be the collection of continuous functions $f : X \rightarrow I$. 

For $f \in C(X,I)$, let $\|f\| = \sup_{x \in X} |f(x)|$ be the uniform norm of $f$.

If $f,g \in C(X,I)$, then define $(f \wedge g)(x) = \min\{f(x),g(x)\}$ and $(f \vee g)(x) = \max\{f(x),g(x)\}$. 

$L \subseteq C(X,I)$ is a lattice if and only if if for all $f,g \in L$, $f \wedge g, f \vee g \in L$. 
\end{definition}

The following fact follows from the proof of the Stone-Weierstass theorem:

\Begin{fact}{stone weierstrass theorem}
Let $X$ be a compact metric space. Let $I \subseteq \bbR$ be an interval, $f \in C(X,I)$, and $L \subseteq C(X,I)$ be a lattice. 

Suppose for all $x,y \in X$, there is some $g \in L$ so that $g(x) = f(x)$ and $g(y) = f(y)$. Then for all $\epsilon > 0$, there is some $h \in L$ so that $\|f - h\| < \epsilon$. 
\end{fact}

\begin{proof}
Fix $f \in C(X,I)$. Pick an $\epsilon > 0$. Fix $x \in X$. By the assumption, for each $y \in X$, choose functions $g_y^x \in L$ so that $f(x) = g_y^x(x)$ and $f(y) = g_y^x(y)$. Define $A^x_y = \{z \in X : g_y^x(z) < f(z) + \epsilon\}$. $A^x_y$ is open using the continuity of $g_y^x$ and $f$. Note that $y \in A^x_y$. $\bigcup_{y \in X} A^x_y$ covers $X$. By compactness, there is a finite set $F_x \subseteq X$ so that $\bigcup_{y \in F_x} A_y^x = X$. 

Define $k_x = \bigwedge_{y \in F_x} g_y^x$. Since $L$ is a lattice, $k_x \in L$. Note that $k_x(x) = f(x)$. For any $z \in X$, there is some $y \in F_x$ so that $z \in A_y^x$. Then $k_x(z) \leq g_x^y(z) < f(z) + \epsilon$. This show that for all $z \in X$, $k_x(z) < f(z) + \epsilon$. 

For $x \in X$, let $B_x = \{z \in X : f(z) - \epsilon < k_x(z)\}$. Each $B_x$ is open and $x \in B_x$. $\bigcup_{x \in X} B_x = X$. By compactness, there is a finite set $F \subseteq X$ so that $\bigcup_{x \in F} B_x = X$. 

Let $h(z) = \bigvee_{x \in F} k_x(z)$. Again since $L$ is a lattice, $h \in L$. Pick any $z \in X$. By the above, $f(z) - \epsilon < k_x(z) \leq h(z)$. Also by the above, $k_x(z) < f(z) + \epsilon$. Hence $\|f - h\| < \epsilon$.
\end{proof}

Fact \ref{uniform convergence of special function} is the main technical approximation that will be needed. The following notation facilitates the exposition.

\Begin{definition}{extension of modulus}
If $\Delta$ is a $k$-ary modulus, then let $\tilde \Delta : \bbR^n \rightarrow \gbbR$ be defined by $\tilde \Delta(\bar{x}) = \Delta(\pi(\bar{x}))$. $\tilde\Delta$ is a continuous function that respects $\Delta$ using subadditivity.

Let $f : \bbR^k \rightarrow \bbR$ be a continuous function and $\Delta$ be a $k$-ary modulus. Let $U^{f,\bar{x},\Delta}(\bar{z}) = f(\bar{x}) + \tilde\Delta(\bar{z} - \bar{x})$. Let $L^{f,\bar{x},\Delta}(\bar{z}) = f(\bar{x}) - \tilde\Delta(\bar{z} - \bar{x})$. Note that $f$ respects $\Delta$ if and only for all $\bar{x}$, $L^{f,\bar{x},\Delta}(\bar{z}) \leq f(\bar{z}) \leq U^{f,\bar{x},\Delta}(\bar{z})$ for all $\bar{z}$. 

Let $\Delta$ be a $k$-ary modulus of uniform continuity so that $\Delta[[0,1]^k] \subseteq [0,1]$. Let $\bar{x},\bar{y} \in [0,1]^k$. Suppose $\Delta(\pi(\bar{y} - \bar{x})) > 0$. Let $a,b \in [0,1]$ with $a \leq b < a + \Delta(\pi(\bar{y} - \bar{x}))$. Let $\tilde\Delta^{\bar{x},\bar{y}}_{a,b} : [0,1]^k \rightarrow [0,1]$ be defined by
$$\tilde\Delta^{\bar{x},\bar{y}}_{a,b}(\bar{z}) = \min\left\{1, a + \frac{b - a}{\tilde\Delta(\bar{y} - \bar{x})}\tilde\Delta(\bar{z} - \bar{x})\right\}$$
Note that since $\tilde\Delta$ respects $\Delta$, so does the $\tilde\Delta_{a,b}^{\bar{x},\bar{y}}$. 

If $\Delta(\pi(\bar{y} - \bar{x})) = 0$, then define
$$\tilde\Delta^{\bar{x},\bar{y}}_{a,a}(\bar{z}) = a.$$
This also respects $\Delta$. 

For $\bar{x}, \bar{y},a,b$ satisfying the above conditions, $\tilde\Delta^{\bar{x},\bar{y}}_{a,b} : [0,1]^k \rightarrow [0,1]$, $\tilde\Delta^{\bar{x},\bar{y}}_{a,b}(\bar{x}) = a$, $\tilde\Delta_{a,b}^{\bar{x},\bar{y}}(\bar{y}) = b$, and this function respects $\Delta$. 

Note that if $\bar{x},\bar{y} \in (\bbQ \cap [0,1])^k$ and $a,b \in \bbQ \cap [0,1]$, then the code for $\tilde\Delta^{\bar{x},\bar{y}}_{a,b}$ is arithmetic in the code for $\Delta$. 
\end{definition}

\Begin{fact}{uniform convergence of special function}
Let $\Delta$ be a $k$-ary-modulus. Let $\bar{x},\bar{y},a,b$ and $\bar{x}',\bar{y}',a',b'$ satisfy the conditions in Definition \ref{extension of modulus}. Let $M$ the maximum of $\Delta$ on $[0,1]^k$. Then
$$\|\tilde\Delta_{a,b}^{\bar{x},\bar{y}} - \tilde\Delta_{a',b'}^{\bar{x}',\bar{y}'}\| \leq |a - a'| + \left|\frac{b - a}{\tilde\Delta(\bar{y} - \bar{x})}\right|\tilde\Delta(\bar{x} - \bar{x}') + M \left|\frac{b - a}{\tilde\Delta(\bar{y} - \bar{x})} - \frac{b' - a'}{\tilde\Delta(\bar{y}' - \bar{x}')}\right|$$
The main observation is that for a fixed $\bar{x},\bar{y},a,b$, the latter expression in the above gets arbitrarily close to $0$ as $\bar{x}',\bar{y}',a',b'$ get close to $\bar{x},\bar{y},a,b$, respectively. Therefore $\|\tilde\Delta_{a,b}^{\bar{x},\bar{y}} - \tilde\Delta_{a',b'}^{\bar{x}',\bar{y}'}\|$ approaches $0$ as $\bar{x}',\bar{y}',a',b'$ gets close to $\bar{x},\bar{y},a,b$.
\end{fact}

\begin{proof}
For any $\bar{z} \in [0,1]^k$, one has by definition
$$|\tilde\Delta_{a,b}^{\bar{x},\bar{y}}(\bar{z}) - \tilde\Delta_{a',b'}^{\bar{x}',\bar{y}'}(\bar{z})| \leq \left|a + \frac{b - a}{\tilde\Delta(\bar{y} - \bar{x})}\tilde\Delta(\bar{z} - \bar{x}) - a' - \frac{b'-a'}{\tilde\Delta(\bar{y}' - \bar{x}')}\tilde\Delta(\bar{z} - \bar{x}')\right|$$
Using the triangle inequality to extract out $|a - a'|$, one obtains
$$\leq |a - a'| + \left|\frac{b - a}{\tilde\Delta(\bar{y} - \bar{x})}\tilde\Delta(\bar{z} - \bar{x}) -\frac{b'-a'}{\tilde\Delta(\bar{y}' - \bar{x}')}\tilde\Delta(\bar{z} - \bar{x}')\right|$$ 
By subtracting and adding the same expression, one has
$$= |a - a'| + \left|\frac{b - a}{\tilde\Delta(\bar{y} - \bar{x})}\tilde\Delta(\bar{z} - \bar{x}) - \frac{b - a}{\tilde\Delta(\bar{y} + \bar{x})}\tilde\Delta(\bar{z}-\bar{x}') + \frac{b - a}{\tilde\Delta(\bar{y} - \bar{x})}\tilde\Delta(\bar{z}-\bar{x}') - \frac{b'-a'}{\tilde\Delta(\bar{y}' - \bar{x}')}\tilde\Delta(\bar{z} - \bar{x}')\right|$$
Using the triangle inequality and factoring, one has
$$= |a - a'| + \left|\frac{b - a}{\tilde\Delta(\bar{y} - \bar{x})}\right||\tilde\Delta(\bar{z} - \bar{x}) - \tilde\Delta(\bar{z} - \bar{x}')| + \tilde\Delta(\bar{z} - \bar{x}') \left|\frac{b - a}{\Delta(\pi(\bar{y} - \bar{x}))} - \frac{b' - a'}{\tilde\Delta(\bar{y}'-\bar{x}')}\right|$$
By assumption, $\tilde\Delta(\bar{z} - \bar{x}') \leq M$. Using the properties of $\Delta$ from Definition \ref{modulus}, one can show that $|\tilde\Delta(\bar{z} - \bar{x}) - \tilde \Delta(\bar{z} - \bar{x}')| = |\Delta(\pi(\bar{z} - \bar{x})) - \Delta(\pi(\bar{z} - \bar{x}'))| \leq |\tilde\Delta(\bar{x} - \bar{x}')|$. Thus one has
$$\leq |a - a'| + \left|\frac{b - a}{\tilde\Delta(\bar{y} - \bar{x})}\right|\tilde\Delta(\bar{x} - \bar{x}') + M \left|\frac{b - a}{\tilde\Delta(\bar{y} - \bar{x})} - \frac{b' - a'}{\tilde\Delta(\bar{y}' - \bar{x}')}\right|$$
The first statement has been verified. 

For the main observation: Note that as $a$ approaches $a'$, the first term goes to $0$. As $\bar{x}$ approaches $\bar{x}'$, the second term approaches $0$. As $a$, $b$, $\bar{x}$, $\bar{y}$ approaches $a'$, $b'$, $\bar{x}'$, $\bar{y}'$, respectively, the third term goes to $0$. 
\end{proof}

\Begin{fact}{sublattice density condition}
Let $A \subseteq A' \subseteq C(X,I)$. Let $L(A)$ and $L(A')$ be the lattice generated by $A$ and $A'$, respectively. Suppose for all $\epsilon > 0$ and $f' \in A'$, there is some $f \in A$ so that $\|f - f'\| < \epsilon$. Then for all $\epsilon > 0$ and $f' \in L(A')$, there is some $f \in L(A)$ so that $\|f - f'\| < \epsilon$. 
\end{fact}

\begin{proof}
Observe that if $\|f - f'\| < \epsilon$ and $\|g - g'\| < \epsilon$, then $\|f \wedge g - f' \wedge g'\| < \epsilon$ and $\|f\vee g - f' \vee g'\| < \epsilon$. The result follows from this observation by induction.
\end{proof}

\Begin{definition}{dense collection in formula respecting modulus}
Let $\Delta$ is a $k$-ary modulus. Let $D_\Delta'$ be the smallest lattice containing $\tilde\Delta^{\bar{x},\bar{y}}_{a,b}$ where $\bar{x},\bar{y} \in (\bbR \cap [0,1])^k$ and $a,b \in \bbR \cap [0,1]$ satisfy the conditions in Definition \ref{extension of modulus} (with respect to $\Delta$).

Let $D_\Delta$ be the smallest lattice containing $\tilde\Delta^{\bar{x},\bar{y}}_{a,b}$ where $\bar{x},\bar{y} \in (\bbQ \cap [0,1])^k$ and $a,b \in \bbQ \cap [0,1]$ satisfying the conditions in Definition \ref{extension of modulus} (with the respect to $\Delta$). $D_\Delta$ is a countable set. $D_\Delta$ has a code which is arithmetic in the code of $\Delta$. 
\end{definition}

\Begin{fact}{density of certain collection}
Let $\Delta$ be a modulus of arity $k$. If $u : [0,1]^k \rightarrow [0,1]$ is a function respecting $\Delta$ and $\epsilon > 0$, then there is a $h \in D_\Delta$ so that $\|u - h\| < \epsilon$. 
\end{fact}

\begin{proof}
Fix $\epsilon > 0$. Fix $\bar{x},\bar{y} \in [0,1]^k$. First suppose that $\Delta(\pi(\bar{y} - \bar{x})) > 0$. Without loss of generality, suppose that $u(\bar{x}) \leq u(\bar{y})$. Since $u$ respects the modulus $\Delta$, $u(\bar{x}) \leq u(\bar{y}) < u(\bar{x}) + \Delta(\pi(\bar{y} - \bar{x}))$. Hence $\tilde\Delta^{\bar{x},\bar{y}}_{u(\bar{x}),u(\bar{y})} \in D_\Delta'$. If $\Delta(\pi(\bar{y} - \bar{x})) = 0$, then $u(\bar{x}) = u(\bar{y})$. In this case, consider $\tilde\Delta^{\bar{x},\bar{y}}_{u(\bar{x}),u(\bar{x})} \in D'_\Delta$. This shows that Fact \ref{stone weierstrass theorem} can be used to find some $h \in D_\Delta'$ so that $\|u - h\| < \epsilon$. 

Using Fact \ref{uniform convergence of special function} and Fact \ref{sublattice density condition}, one can find some $h \in D_\Delta$ so that $\|u - h\| < \epsilon$. 
\end{proof}

\Begin{fact}{dense basic formula respecting modulus}
Let $\Omega$ be a weak modulus. Let $\varphi_0, ..., \varphi_{k - 1}$ be atomic formulas of $L_{\omega,\omega}$ with free variables among $v_0, ..., v_{n - 1}$. Then there is a countable set $E(\varphi_0, ..., \varphi_{k - 1})$ with the following properties:

\noindent 1. All formulas $\psi \in E(\varphi_0, ..., \varphi_{k - 1})$ are basic $n$-ary respecting $\Omega$. 

\noindent 2. For any $\epsilon > 0$, any $\SCRL$-structure $\CM$, and any formula of the form $u(\varphi_0, ..., \varphi_{k - 1})$ respecting $\Omega$ where $u : [0,1]^k \rightarrow [0,1]$ is continuous, there is some formula $\psi \in E(\varphi_0, ..., \varphi_{k -1})$ so that $\|\psi^\CM - u(\varphi_0, ..., \varphi_{k - 1})^\CM\| < \epsilon$. 

\noindent 3. There is a real arithmetic in $\hat\SCRL\oplus\hat\Omega$ which codes $E(\varphi_0, ..., \varphi_{k - 1})$. 
\end{fact}

\begin{proof}
Let $\Delta_i$ be the canonical modulus for $\varphi_i$. Suppose some of the $\Delta_i$ are constantly $0$. Without loss of generality, assume that there is some $j < k$ so that $\Delta_i$ is constantly $0$ for all $i \geq j$. Fix a $\SCRL$-structure $\CM$. There are constants $b_j, ..., b_{k - 1}$ in $[0,1]$ so that for each $i \geq j$, $\varphi^\CM_i$ takes constant value $b_j$. Define $\tilde u(\varphi_0, ..., \varphi_{j - 1}) = u(\varphi_0, ..., \varphi_{j - 1}, b_j, ..., b_{k - 1})$. $\tilde u(\varphi_0, ..., \varphi_{j - 1})$ respects the weak modulus $\Omega$ and, $\tilde u(\varphi_0, ..., \varphi_{j - 1})$ and $u(\varphi_0, ..., \varphi_{k - 1})$ are equal when interpreted in $\CM$.

From the above discussion, it suffices to consider the case when all the modulus $\Delta_0, ..., \Delta_{k - 1}$ are not constantly zero.

Let $\Delta$ be the modulus from Fact \ref{largest modulus to respect a modulus}. Let $E(\varphi_0, ..., \varphi_{k - 1})$ be the collection of the formulas of the form $u(\varphi_0, ..., \varphi_{k - 1})$ where $u \in D_\Delta$. 
\end{proof}

Finally, the next fact shows that there is a countable dense set $F_\Omega^n$ of $n$-ary basic formulas respecting the weak modulus $\Omega$ which is dense among the collection of all basic $n$-ary formulas $\varphi$. Also $F_\Omega^n$ has a code arithmetic in $\hat{\SCRL}\oplus \hat{\Omega}$.

\Begin{fact}{dense set of Omega respecting basic formulas}
There is a countable set $F^n_\Omega$ of $n$-ary basic formulas respecting the weak modulus $\Omega$ such that for all $\epsilon > 0$, $\SCRL$-structures $\CM$, and basic $n$-ary-formulas $\varphi$ respecting $\Omega$, there is some $\psi \in F^n_\Omega$ so that $\|\psi^\CM - \varphi^\CM\| < \epsilon$. Moreover, the code for $F^n_\Omega$ is arithmetic in $\hat\SCRL\oplus\hat\Omega$. 
\end{fact}

\begin{proof}
Fix an $\SCRL$-recursive enumeration of all finite tuples $(\varphi_0, ..., \varphi_{k - 1})$ (where $k \in \omega$) of atomic formulas in the free variables $v_0, ..., v_{n - 1}$. Let $F^n_\Omega$ be the union of all such $E(\varphi_0, ..., \varphi_{k - 1})$. This works using Fact \ref{standard form of basic formula} and Fact \ref{dense basic formula respecting modulus}.
\end{proof}

\section{Continuous Scott Analysis}\label{Continuous Scott Analysis}

The following is the back-and-forth pseudo-distance:

\Begin{definition}{back-and-forth pseudo-distance}
(\cite{Metric-Scott-Analysis} Definition 3.1) Let $\SCRL$ be a language of continuous logic. Let $\CM$ be a $\SCRL$-structure. Let $\bar{a},\bar{b}$ be tuples from $M$ of the same length. Let $\Omega$ be a weak modulus. 

Let 
$$r_0(\bar{a},\bar{b}) = \sup_{\varphi}|\varphi^\CM(\bar{a}) - \varphi^\CM(\bar{b})|$$
where the supremum is taken over all basic $\SCRL$-formulas respecting $\Omega$. 

Suppose $r_\beta$ has been defined for all $\beta < \alpha$ and $\alpha$ is a limit ordinal, then let
$$r_\alpha(\bar{a},\bar{b}) = \sup_{\beta < \alpha} r_\beta(\bar a,\bar b)$$

Suppose $r_\alpha$ has been defined, let
$$r_{\alpha + 1}(\bar{a},\bar{b}) = \sup_{c,d \in M}\inf_{c',d' \in M} r_\alpha(\bar{a}c, \bar{b}d') \vee r_\alpha(\bar{a}c',\bar{b}d')$$

For each $n \in \omega$, let $r_{\alpha,n}$ be the restriction of $r_\alpha$ to $M^n \times M^n$. 
\end{definition}

\Begin{fact}{continuity aspect of pseudo-distance}
(\cite{Metric-Scott-Analysis} Lemma 3.2 and 3.3) Let $\alpha$ be an ordinal and $n \in \omega$. $r_{\alpha,n}$ is a pseudo-distance. Each $r_{\alpha,n}$ is a uniformly continuous function respecting the modulus $\Omega_n$. 

If $\alpha < \beta$ and $\bar{a}$ is a tuple in $M$, then $r_\alpha(\bar a,\bar b) \leq r_\beta(\bar a, \bar b)$.
\end{fact}

\Begin{definition}{continuous scott rank}
(\cite{Metric-Scott-Analysis} Definition 3.6) The least ordinal $\alpha$ so that $r_\alpha = r_{\alpha + 1}$ is the Scott rank of $\CM$ with respect to $\Omega$ and is denoted $\mathrm{SR}_\Omega(\CM)$. 
\end{definition}

\Begin{fact}{computation of r0 in admissible set}
Let $\SCRL$ be a countable language with code $\hat\SCRL$. Let $\CD$ be a countable $\SCRL$-pre-structure. Let $\bar{D}$ denotes its completion structure. Let $\Omega$ be a weak modulus with code $\hat\Omega$. Let $\CA$ be a countable admissible set containing $\hat\SCRL\oplus\hat\Omega\oplus\CD$. Let $\bar{a},\bar{b}$ be tuples from $D$ of the same length. Let $(r_0(\bar{a},\bar{b}))^\CA$ be the computation of $r_0(\bar{a},\bar{b})$ in $\CA$. Then $r_0(\bar{a},\bar{b}) = (r_0(\bar{a},\bar{b}))^\CA$. 
\end{fact}

\begin{proof}
Fact \ref{admissible has correct interpretation} implies that for any $\varphi \in L_{\omega,\omega}^\CA$, $\varphi^{\bar\CD}(\bar{a}) = (\varphi^{\bar{\CD}}(\bar{a}))^\CA$ and $\varphi^{\bar\CD}(\bar{b}) = (\varphi^{\bar\CD}(\bar{b}))^\CA$. That is, the computation of formulas $\varphi \in L_{\omega,\omega}^\CA$ is the same as in the real world.

Let $n$ be the length $\bar{a}$. The computation of $r_0(\bar{a},\bar{b})$ in $\CA$ entails using only the basic $n$-ary formulas respecting $\Omega$ which belong to $\SCRL_{\omega,\omega}^\CA$ when evaluating the supremum. The set $F_{\Omega}^n$ from Fact \ref{dense set of Omega respecting basic formulas} is contained in and is dense in the collection of all $n$-ary basic formulas respecting $\Omega$ under the uniform norm. The entire set $F_{\Omega}^n$ is hyperarithmetic in $\hat\SCRL\oplus\hat\Omega$. So $F_{\Omega}^n \in \CA$. In particular, each function of $F_{\Omega}^n$ is in $L_{\omega,\omega}^\CA$. 

Together, these facts show that the computation of $r_0(\bar{a},\bar{b})$ in $\CA$ gives the true value of $r_0(\bar{a},\bar{b})$. 
\end{proof}

\Begin{theorem}{bounds on continuous scott rank}
Let $\SCRL$ be a countable language with code $\hat\SCRL$. Let $\Omega$ be a weak modulus with code $\hat\Omega$. Let $\CD$ be a countable $\SCRL$-pre-structure. Let $\bar{\CD}$ denote its completion structure. Then $\mathrm{SR}_\Omega(\bar\CD) \leq \omega_1^{\hat\SCRL\oplus\hat\Omega\oplus\CD}$. 
\end{theorem}

\begin{proof}
Let $\CA$ be any admissible set containing $\hat\SCRL \oplus \hat\Omega \oplus \CD$. 

Fix a $q \in \bbQ^{> 0}$. Define the following operator $\Gamma_q : \mathscr{P}(D^{<\omega} \times D^{<\omega}) \rightarrow \mathscr{P}(D^{<\omega} \times D^{<\omega})$ by: $\Gamma(X)$ is the set of $(\bar{a},\bar{b})$ so that one of the following hold:

1. $|\bar{a}| \neq |\bar{b}|$

2. $r_0(\bar{a},\bar{b}) > q$.

3. $(\exists c \in D)(\exists d \in D)(\forall c' \in D)(\forall d' \in D)((\bar{a}c',\bar{b}d) \in X \vee (\bar{a}c,\bar{b},d') \in X)$

For each $q \in \bbQ$, $\Gamma_q$ is defined by a positive $\Sigma$-formula with parameters from $\CA$. By results of Gandy (see \cite{Admissible-Sets-and-Structures}, Chapter VI, Corollary 2.8), the closure ordinal $\|\Gamma_q\|$ of $\Gamma_q$ is less than or equal to the $o(\CA) = A \cap \mathrm{ON}$, the ordinal height of $\CA$. 

Let $I_{\Gamma_q}^\infty = I_{\Gamma_q}^{\|\Gamma_q\|}$ be the least fixed point of $\Gamma_q$. By Fact \ref{computation of r0 in admissible set}, $r_0$ is computed correctly by $\CA$. So $(\bar{a},\bar{b}) \in I_{\Gamma_q}^\infty$ if and only if there is some $\alpha$ so that $r_\alpha(\bar{a},\bar{b}) > q$. 

Let $\gamma = \sup_{q \in \bbQ^{> 0}} \|\Gamma_q\|$. By the above, $\gamma \leq o(\CA)$. So $r_\gamma(\bar{a}, \bar{b}) = r_\beta(\bar{a}, \bar{b})$ for all tuples $\bar{a}, \bar{b}$ from $D$ and $\beta \geq \gamma$. By Fact \ref{continuity aspect of pseudo-distance}, each $r_{\beta,n}$ is a continuous function on $\bar{D}^n \times \bar{D}^n$. Since for each $\beta \geq \gamma$, $r_{\gamma,n}\upharpoonright D^n \times D^n = r_{\beta,n}\upharpoonright D^n \times D^n$, one has that $r_{\gamma,n} \upharpoonright \bar{D}^n \times \bar{D}^n = r_{\beta,n} \upharpoonright \bar{D}^n \times \bar{D}^n$ for all $\beta \geq \gamma$. Hence $r_\beta(\bar{a}, \bar{b}) = r_\gamma(\bar{a}, \bar{b})$ for all $\beta \geq \gamma$ and all tuples $\bar{a}, \bar{b}$ of the same length from $\bar{D}$. $\mathrm{SR}_\Omega(\bar{\CD}) \leq \gamma$ and in fact is equal. As $\gamma \leq o(\CA)$, one has $\mathrm{SR}_\Omega(\bar{\CD}) \leq o(\CA)$. 

Letting $\CA = L_{\omega_1^{\hat\SCRL\oplus\hat\Omega\oplus\CD}}(\hat\SCRL\oplus\hat\Omega\oplus\CD)$, which has ordinal height $o(\CA) = \omega_1^{\hat\SCRL\oplus\hat\Omega\oplus\CD}$, completes the proof.
\end{proof}

\bibliographystyle{amsplain}
\bibliography{references}

\end{document}